\documentclass[12pt,twoside]{article}

\date{}

\begin{document}

\centerline{}

\centerline {\Large{\bf A Note on Intuitionistic Fuzzy Hypervector Spaces}}

\newcommand{\mvec}[1]{\mbox{\bfseries\itshape #1}}
\centerline{}
\centerline{\bf {Sanjay Roy, T. K. Samanta }}

\centerline{}

\centerline{Department of Mathematics, South Bantra Ramkrishna Institution , West Bengal, India. }
\centerline{e-mail: sanjaypuremath@gmail.com}

\centerline{Department of Mathematics, Uluberia
College, West Bengal, India.}
\centerline{e-mail: mumpu$_{-}$tapas5@yahoo.co.in}

\newtheorem{Theorem}{\quad Theorem}[section]

\newtheorem{definition}[Theorem]{\quad Definition}

\newtheorem{theorem}[Theorem]{\quad Theorem}

\newtheorem{remark}[Theorem]{\quad Remark}

\newtheorem{corollary}[Theorem]{\quad Corollary}

\newtheorem{note}[Theorem]{\quad Note}

\newtheorem{lemma}[Theorem]{\quad Lemma}

\newtheorem{example}[Theorem]{\quad Example}

\newtheorem{result}[Theorem]{\quad Result}

\newtheorem{proposition}[Theorem]{\quad Proposition}

\begin{abstract}
\textbf{\emph{ The notion of Intuitionistic fuzzy hypervector space has been generalized and a few basic properties on this concept are studied. It has been shown that the intersection and union of an arbitrary family of Intuitionistic fuzzy hypervector spaces are also Intuitionistic fuzzy hypervector space. Lastly, the notion of a linear transformation on a hypervector space is introduced and established an important theorem relative to Intuitionistic fuzzy hypervector spaces.}}
\end{abstract}


\section{Introduction}
\smallskip\hspace{.5cm}The notion of hyperstructure was introduced by F. Marty in 1934. Then he established the definition of hypergroup \cite{Marty1}  in 1935. Since then many researchers have studied and developed ( for example see \cite{Marty2} , \cite{Nakassis} ) the concept of different types of hyperstructures in different views. In 1990 M. S. Tallini \cite{Tallini} introduced the notion of  hypervector spaces. Then in 2005 R. Ameri \cite{Ameri} also studied this spaces extensively. In our previous papers ( \cite{Roy1}, \cite{Roy2} ), we also introduced the notion of a hypervector spaces in more general form than the previous concept of hypervector space and thereafter established a few useful theorems in this space.\\
\smallskip\hspace{.5cm}The concept of intuitionistic Fuzzy set, as a generalization of a fuzzy set was first introduced by Atanassov \cite{Atanassov}. Then many researchers (\cite{dinda}, \cite{Samanta}, \cite{Vijayabalaji}) applied this notion to norm, Continuity and Uniform Convergence etc. At the present time  many researchers ( for example \cite{Torkzadeh} ) are trying to apply this concept on the hyperstructure theory.\\
\smallskip\hspace{.5cm}In this paper, the concept of Intuitionistic fuzzy hypervector space is introduced and a few basic properties are developed. Further it has been shown that the intersection and union of a arbitrary family of Intuitionistic fuzzy hypervector spaces are also Intuitionistic fuzzy hypervector space. Lastly we have introduced the notion of a linear transformation on a hypervector space and established an important theorem relative to Intuitionistic fuzzy hypervector space.

\section{Preliminaries}
\smallskip\hspace{.5cm}This section contain some basic definition and preliminary results which will be needed.

\begin{definition} \cite{Nakassis}
A \textbf{hyperoperation} over a non empty set X is a mapping of \,X$\times $X \,into the set of all non empty subsets of X.
\end{definition}

\begin{definition}\cite{Nakassis}
A non empty set R with exactly one hyperoperation $'\#'$ is a \textbf{hypergroupoid}.\\
\smallskip\hspace{.5cm}Let $(X\, ,\, \#)$ be a hypergroupoid. For every point x $\in$ X and every non empty subset A of X , we defined  $x\;\#\; A\;=\;\bigcup_{a\,\in\,A}\{ x \;\#\; a\}$.
\end{definition}
\begin{definition}\cite{Nakassis}
A hypergroupoid $(X \,,\,\#)$ is called a \textbf{hypergroup} if\\
$( i )$\hspace{0.5cm}   $x\;\#\;(y \;\#\;z)=(x\;\#\;y)\;\#$ z\\
$( ii )$\hspace{0.4cm}  $\exists\, 0\,\in$ X such that for every a $\in$ X , there is unique element b $\in$ X for which $0$ $\in$ a $\#$ b and $0\,\in$ b $\#$ a . Here b is denoted by $-$a .\\
$( iii )$ \hspace{0.3cm}For all a , b , c $\in$ X if a $\in$ b $\#$ c , then b $\in\; a\; \#\; (- c)$.
\end{definition}

\begin{proposition}\cite{Nakassis}
$(i)$ In a hypergroup $(X \,,\, \#)$ , $-(- a) = a$ , $\forall  a \,\in$ X.

$(ii)$ $0 \;\#\; a\,=\,\{a\}$ , $\forall a \in$ X , if $(X\,,\,\#)$ is a commutative hypergroup.

$(iii)$ In a commutative hypergroup $($X , $\#)$ , $0$ is unique.
\end{proposition}
\begin{definition}\cite{Roy1}
A \textbf{hyperring} is a non empty set equipped with a hyperaddition '$\#$' and a multiplication `.' such that $(X\,,\,\#)$ is a commutative hypergroup and $(X\,,.)$ is a semigroup and the multiplication is distributive across the hyperaddition both from the left and from the right and \,$a.0 = 0.a = 0, \,\forall \,a\in$ X , where $0$ is the zero element of the hyperring.
\end{definition}
\begin{definition}\cite{Roy1}
A \textbf{hyperfield} is a non empty set X equipped with a hyperaddition `\,$\#$' and a multiplication `.'such that\\
$(\,i\,)$\hspace{0.5cm}$(X\,,\;\#\;,\,.)$ is a hyperring.\\
$(\,ii\,)$\hspace{0.4cm}$\exists$ an element 1 $\in$ X, called the identity element such that a.1 = a, \\$\forall \,a\in$ X \\
$(\,iii\,)$\hspace{0.3cm}For each non zero element a in X , $\exists$ an element $a^{-1}$ such that $a.a^{-1}$=1\\
$(\,iv\,)$\hspace{0.4cm}a.b = b.a , $\forall$ a , b\,$\in$X.
\end{definition}

\begin{definition}\cite{Roy2}
Let $(F\,,\oplus\,,.)$ be a hyperfield and $(V\,,\,\#)$ be an additive commutative hypergroup . Then $V$ is said to be a \textbf{ hypervector space} over the hyperfield $F$ if there exist a hyperoperation $\ast: \;F \;\times\, V \,\rightarrow\;P^{\ast}\,(\,V\,)$ such that\\
$(\,i\,)$\hspace{0.5cm}a $\ast\,(\alpha\;\#\;\beta)\,\subseteq \,a\ast\alpha\;\,\#\;\,a\ast\beta$ , \hspace{1cm}$\forall\, a\in F\;and\;\forall\, \alpha\, , \beta \in$ $V$\\
$(\,ii\,)$\hspace{0.4cm}$(a\,\oplus\,b)\ast\alpha \subseteq a\ast\alpha\;\#\;b\ast\alpha$  , \hspace{1cm} $\forall\,a , b\in F\;and\;\forall\,\alpha\in
$ $V$\\
$(\,iii\,)$\hspace{0.3cm}$(a\,.\,b)\ast\alpha\,=\,a\ast(b\ast\alpha)$  , \hspace{2cm} $\forall\, a , b \in$F and $\forall\, \alpha\in$ $V$.\\
$(\,iv\,)$\hspace{0.4cm}$(-a)\ast\alpha\,=\,a\ast(-\alpha) ,\hspace{2.7cm}$  $\forall\, \alpha\in$ V and $\forall\, a\in$ $F$.\\
$(\,v\,)$\hspace{0.5cm}$\alpha\in 1_{F}\ast\alpha$,  $\theta\in\, 0 \ast\alpha$ and $0\ast\theta\,=\,\theta,$ \hspace{1cm}$\forall\, \alpha\in$ $V$.\\
 where $1_{F}$ is the identity element of $F$ , $0$ is the zero element of $F$ and $\theta$ is zero vector of $V$ and $P^{\ast}(\,V\,)$ is the set of all non empty subset of $V$.
 \end{definition}

\begin{definition}\cite{Ameri}
Let f : X $\rightarrow$ Y be a mapping and $\nu\in$ FS(Y). Then we define $f^{-1}(\nu)\in$ FS(X) as follows:\\
$f^{-1}(\nu)(x)\,=\,\nu(f(x))$, $\forall \;x\in$ X.
\end{definition}
\begin{definition}\cite{Vijayabalaji}
Let E be a any set. An Intuitionistic fuzzy set ( IFS ) A of E is an object of the form A$\,=\{\,(\,x,\,\mu_{A}(x),\,\nu_{A}(x)\,)\;:\;\;x\,\in\,E\,\}$, where the functions $\mu_{A}\,:\,E\,\rightarrow\,[0\,,\,1]$ and $\nu_{A}\,:\,E\,\rightarrow\,[0\,,\,1]$ denotes the degree of membership and the non-membership of the element $x\,\in\,E$ respectively and for every $x\,\in\,E$, $0\,\leq\,\mu_{A}(x)\,+\,\nu_{A}(x)\,\leq\,1$.
\end{definition}

\section{Intuitionistic Fuzzy hypervector Space}
\smallskip\hspace{.5cm}In this section we established the definition of intuitionistic fuzzy hypervector Spaces and deduce some important theorems.

\begin{definition}
Let $(F,\oplus,.)$ be a hyperfield. An \textbf{intuitionistic fuzzy hyperfield} on $F$ is an object of the form  $\;A\,=\{\,(\,a,\,\mu_{F}(a),\,\nu_{F}(a)\,)\;:\;\;a\,\in\,F\,\}\,$ satisfies the following conditions\,:\\
$(i)$ $\bigwedge_{x\in a\oplus b}$ $\mu_{F}(x)\;\geq\, \mu_{F}(a)\wedge\mu_{F}(b)$, $\forall$ a, b\,$\in$ F\\
$(ii)$ $\mu_{F}(-a)\,\geq \,\mu_{F}(a)$, $\forall\; a\,\in$ F\\
$(iii)$ $\mu_{F}(a.b)\,\geq\, \mu_{F}(a)\wedge\mu_{F}(b)$ $\forall$ a, b $\in$ F\\
$(iv)$ $\mu_{F}(a^{-1})\,\geq \,\mu_{F}(a)$, $\forall \;a(\neq 0)\,\in$ F\\
$(v)$ $\bigvee_{x\in a\oplus b}$ $\nu_{F}(x)\;\leq\, \nu_{F}(a)\vee\nu_{F}(b)$, $\forall$ a, b $\in$ F\\
$(vi)$ $\nu_{F}(-a)\,\leq \,\nu_{F}(a)$, $\forall\; a\,\in$ F\\
$(vii)$ $\nu_{F}(a.b)\,\leq\,\nu_{F}(a)\vee\nu_{F}(b)$, $\forall$ a, b $\in$ F\\
$(viii)$ $\nu_{F}(a^{-1})\,\leq\, \nu_{F}(a)$, $\forall \;a(\neq 0)\,\in$ F.\\

\end{definition}

\begin{result}
If A is a intuitionistic fuzzy hyperfield of F, then\\
$(i)$ $\mu_{F}(0)\,\geq\,\mu_{F}(a)$, $\forall$ a $\in$ F\\
$(ii)$ $\mu_{F}(1)\,\geq\,\mu_{F}(a)$, $\forall$ a $\in\, F\setminus \{0\}$\\
$(iii)$ $\mu_{F}(0)\,\geq\,\mu_{F}(1)$\\
$(iv)$ $\nu_{F}(0)\,\leq\, \nu_{F}(a)$, $\forall$ a $\in$ F\\
$(v)$ $\nu_{F}(1)\,\leq\,\nu_{F}(a)$, $\forall$ a $\in F\setminus \{0\}$\\
$(vi)$ $\nu_{F}(0)\,\leq \,\nu_{F}(1)$
\end{result}

\textbf{Proof} : \,Obvious.

\begin{definition}
Let $(V,\#,\ast)$ be a hypervector space over a hyperfield $(F,\oplus,.)$ and A be a intuitionistic fuzzy hyperfield in F. A  intuitionistic fuzzy subset B$\,=\{\,(\,x,\,\mu_{V}(x),\,\nu_{V}(x)\,)\;:\;\;x\,\in\,V\,\}$ of V is said to be a \textbf{intuitionistic fuzzy hypervector space} of V over a intuitionistic fuzzy hyperfield A, if the following conditions are satisfied:\\
$(i)$ $\bigwedge_{\alpha\in x\# y}$ $\mu_{V}(\alpha)\;\geq \mu_{V}(x)\wedge\mu_{V}(y)$, $\forall\, x,\, y \,\in$ V\\
$(ii)$ $\mu_{V}(-x)\,\geq \,\mu_{V}(x)$, $\forall\; x\;\in$ V\\
$(iii)$ $\bigwedge_{y\in a\ast x}\mu_{V}(y)\,\geq\, \mu_{V}(x)\wedge\mu_{F}(a)$, $\;\forall$ a $\in$ F and $\forall\, x \,\in$ V\\
$(iv)$ $\mu_{F}(1)\,\geq\, \mu_{V}(\theta)$, where $\theta$ be the null vector of V.\\
$(v)$ $\bigvee_{\alpha\in x\# y}$ $\nu_{V}(\alpha)\;\leq \,\nu_{V}(x)\vee\nu_{V}(y)$, $\forall\, x,\, y \,\in$ V\\
$(vi)$ $\nu_{V}(-x)\,\leq\, \nu_{V}(x)$, $\forall\; x\;\in$ V\\
$(vii)$ $\bigvee_{y\in a\ast x}\nu_{V}(y)\,\leq\, \nu_{V}(x)\vee \nu_{F}(a)$, $\forall$ a $\in$ F and $\forall\, x \,\in$ V\\
$(viii)$ $\nu_{F}(1)\,\leq\, \nu_{V}(\theta)$, where $\theta$ be the null vector of V.\\

Here we say that B is a intuitionistic fuzzy hypervector space over a intuitionistic fuzzy hyperfield A.
\end{definition}

\begin{result}
If B is a intuitionistic fuzzy hypervector space over a intuitionistic fuzzy hyperfield A, then\\
$(i)$ $\mu_{F}(0)\,\geq\,\mu_{V}(\theta)$, \\
$(ii)$ $\mu_{V}(\theta)\,\geq\,\mu_{V}(x)$, $\forall\, x \,\in$ V\\
$(iii)$ $\mu_{F}(0)\,\geq\,\mu_{V}(x)$, $\forall\, x \,\in$ V\\
$(iv)$ $\nu_{F}(0)\,\leq \,\nu_{V}(\theta)$, \\
$(v)$ $\nu_{V}(\theta)\,\leq\, \nu_{V}(x)$, $\forall\, x \,\in$ V\\
$(vi)$ $\nu_{F}(0)\,\leq\, \nu_{V}(x)$, $\forall\,x \,\in$ V
\end{result}

\textbf{Proof} : Obvious.

\begin{theorem}
Let V be a hypervector space over a hyperfield F and A be a intuitionistic fuzzy hyperfield. Let B $\in$ IFS(V). Then B is a intuitionistic fuzzy hypervector space over A iff\\
$(i)$ $\bigwedge_{z\in a\ast x\# b\ast y} \mu_{V}(z)\geq (\mu_{F}(a)\wedge\mu_{V}(x))\wedge(\mu_{F}(b)\wedge\mu_{V}(y))$ , $\forall\, x, \,y \,\in$ V and \\$\forall$ $a, \,b\,\in$ F\\
$(ii)$ $\mu_{F}(1)\geq\mu_{V}(\theta)$ where $\theta$ be the null vector of V.\\
$(iii)$ $\bigvee_{z\in a\ast x\# b\ast y} \nu_{V}(z)\leq (\nu_{F}(a)\vee\nu_{V}(x))\vee(\nu_{F}(b)\vee\nu_{V}(y))$ , $\forall\, x,\, y \,\in$ V and \\$\forall$ $a,\, b\,\in$ F\\
$(iv)$ $\nu_{F}(1)\leq \nu_{V}(\theta)$ where $\theta$ be the null vector of V.
\end{theorem}

\textbf{Proof}: First we suppose that B is a intuitionistic fuzzy hypervector space over the intuitionistic fuzzy hyperfield A. Then for a, b $\in$ F and $x,\, y \,\in$ V, we have\\
$\bigwedge_{z\in a\ast x\# b\ast y} \mu_{V}(z)$ =  $\bigwedge_{z\in\alpha\#\beta,\alpha\in a\ast x,\beta\in b\ast y}\mu_{V}(z)$\\
\smallskip$\hspace{3cm}$=$\bigwedge_{\alpha\in a\ast x,\beta\in b\ast y}(\bigwedge_{z\in\alpha\#\beta}\mu_{V}(z))$\\
\smallskip$\hspace{3cm}\geq\bigwedge_{\alpha\in a\ast x,\beta\in b\ast y}(\mu_{V}(\alpha)\wedge\mu_{V}(\beta))$\\
\smallskip$\hspace{3cm}=(\;\;\bigwedge_{\alpha\in a\ast x} \mu_{V}(\alpha))\bigwedge(\;\;\bigwedge_{\beta\in b\ast y} \mu_{V}(\beta))$\\
\smallskip$\hspace{3cm}\geq (\mu_{F}(a)\wedge\mu_{V}(x))\bigwedge(\mu_{F}(b)\wedge\mu_{V}(y))$\\
$\bigvee_{z\in a\ast x\# b\ast y} \nu_{V}(z)$ =  $\bigvee_{z\in\alpha\#\beta,\alpha\in a\ast x,\beta\in b\ast y}\nu_{V}(z)$\\
\smallskip$\hspace{3cm}$=$\bigvee_{\alpha\in a\ast x,\beta\in b\ast y}(\bigvee_{z\in\alpha\#\beta}\nu_{V}(z))$\\
\smallskip$\hspace{3cm}\leq \bigvee_{\alpha\in a\ast x,\beta\in b\ast y}(\nu_{V}(\alpha)\vee\nu_{V}(\beta))$\\
\smallskip$\hspace{3cm}=(\;\;\bigvee_{\alpha\in a\ast x} \nu_{V}(\alpha))\bigvee(\;\;\bigvee_{\beta\in b\ast y} \nu_{V}(\beta))$\\
\smallskip$\hspace{3cm}\leq (\nu_{F}(a)\vee\nu_{V}(x))\bigvee(\nu_{F}(b)\vee\nu_{V}(y))$\\
The second and fourth inequalities are directly follow.
\\
\textbf{conversely} suppose that the inequalities of the theorem hold for all $x,\, y\,\in$ V and $\forall$ a, b $\in$ F.\\
Then $\bigwedge_{z\in 1\ast x\#1\ast y}\mu_{V}(z)\,\geq \,(\mu_{F}(1)\wedge\mu_{V}(x))\bigwedge(\mu_{F}(1)\wedge\mu_{V}(y))$\\
\smallskip$\hspace{4.2cm}\geq \,(\mu_{V}(\theta)\wedge\mu_{V}(x))\bigwedge(\mu_{V}(\theta)\wedge\mu_{V}(y))$\\
\smallskip$\hspace{4.2cm}=\,\mu_{V}(x)\wedge\mu_{V}(y)$\\
i.e $\bigwedge_{z\in  x\# y}\mu_{V}(z)\,\geq\,\mu_{V}(x)\wedge\mu_{V}(y)$, as $x\,\in\,1\ast x$ and $y\,\in\,1\ast y$\\
$\mu_{V}(-x)\,=\,\bigwedge_{z\in -1\ast x}\mu_{V}(z)$  , as $-x\,\in\,-1\ast x$\\
\smallskip$\hspace{1.5cm}\geq \,\bigwedge_{z\in -1\ast x\#0\ast x}\mu_{V}(z)\geq (\mu_{F}(-1)\wedge\mu_{V}(x))\bigwedge(\mu_{F}(0)\wedge\mu_{V}(x))$\\
\smallskip$\hspace{5cm}\geq\, (\mu_{F}(1)\wedge\mu_{V}(x))\bigwedge\mu_{V}(x)$\\
\smallskip$\hspace{5cm}\geq\, (\mu_{V}(\theta)\wedge\mu_{V}(x))\bigwedge\mu_{V}(x)$\\
\smallskip$\hspace{5cm}=\,\mu_{V}(x)\wedge\mu_{V}(x)$\\
\smallskip$\hspace{5cm}=\,\mu_{V}(x)$\\
i.e $\mu_{V}(-x)\geq\mu_(x)$\\
$\bigwedge_{y\in a\ast x}\mu_{V}(y)\,\geq\,\bigwedge_{y\in a\ast x\# 0\ast x}\mu_{V}(y)\,\geq\,(\mu_{F}(a)\wedge\mu_{V}(x))\bigwedge(\mu_{F}(0)\wedge\mu_{V}(x))$\\
\smallskip$\hspace{6cm}=\,(\mu_{V}(x)\wedge\mu_{F}(a))\bigwedge\mu_{V}(x)$  , as $\mu_{F}(0)\geq\mu_{V}(x)$ \\
\smallskip$\hspace{6cm}=\,\mu_{V}(x)\wedge\mu_{F}(a)$\\
The fourth inequality of definition 3.3 is directly follows.\\
Next $\bigvee_{z\in 1\ast x\#1\ast y}\nu_{V}(z)\,\leq\, (\nu_{F}(1)\vee\nu_{V}(x))\bigvee(\nu_{F}(1)\vee\nu_{V}(y))$\\
\smallskip$\hspace{4.2cm}=\, (\nu_{V}(\theta)\vee\nu_{V}(x))\bigvee(\nu_{V}(\theta)\vee\nu_{V}(y))$\\
\smallskip$\hspace{4.2cm}=\,\nu_{V}(x)\vee\nu_{V}(y)$\\
i.e $\bigvee_{z\in  x\# y}\nu_{V}(z)\,\leq\,\nu_{V}(x)\vee\nu_{V}(y)$, as $x\,\in\,1\ast x$ and $y\,\in\,1\ast y$\\
$\nu_{V}(-x)\,\leq\,\bigvee_{z\in -1\ast x}\nu_{V}(z)$  , as $-x\,\in\,-1\ast x$\\
\smallskip\hspace{1.5cm}$\leq\,\bigvee_{z\in -1\ast x\,\# \,0\ast x}\nu_{V}(z)$ \\
\smallskip\hspace{1.5cm}$\leq \,(\nu_{F}(-1)\vee\nu_{V}(x))\vee(\nu_{F}(0)\vee\nu_{V}(x))$\\
\smallskip\hspace{1.5cm}$\leq \,(\nu_{F}(1)\vee\nu_{V}(x))\bigvee \nu_{V}(x)$\\
\smallskip\hspace{1.5cm}$\leq\, (\nu_{V}(\theta)\vee\nu_{V}(x))\bigvee\nu_{V}(x)$\\
\smallskip\hspace{1.5cm}$=\,\nu_{V}(x)\vee\nu_{V}(x)$\\
\smallskip\hspace{1.5cm}$=\,\nu_{V}(x)$\\
i.e $\nu_{V}(-x)\,\leq\,\nu_(x)$\\
$\bigvee_{y\in a\ast x}\nu_{V}(y)\,\leq\, \bigvee_{y\in a\ast x\#\,0\ast x}\nu_{V}(y)$\\
\smallskip\hspace{2.5cm}$\leq\,(\nu_{F}(a)\vee\nu_{V}(x))\vee(\nu_{F}(0)\vee\nu_{V}(x))$\\
\smallskip\hspace{2.5cm}= $(\nu_{V}(x)\vee\nu_{F}(a))\bigvee\nu_{V}(x)$  , as $\nu_{F}(0)\leq\nu_{V}(x)$ \\
\smallskip\hspace{2.5cm}= $\nu_{V}(x)\vee\nu_{F}(a)$\\
The eighth inequality  of definition 3.3 is obvious.\\
Therefore B is a intuitionistic fuzzy hypervector space over A.\\
This completes the proof.
\begin{definition}
Let $B^{\alpha}_{(\alpha\in\Lambda)}\,=\,\{\,x,\,\mu_{V}^{\alpha}(x),\,\nu_{V}^{\alpha}(x)\;:\;\;x\,\in\;V\}$ be a family of intuitionistic fuzzy hypervector spaces of a hypervector space V over the  same intuitionistic fuzzy hyperfield $A\,=\{\,(\,x,\,\mu_{F}(x),\,\nu_{F}(x)\,)\;:\;\;x\,\in\,F\,\}$.Then \\
the \textbf{intersection} of those intuitionistic fuzzy hypervector spaces is defined as\\
$(\bigcap_{\alpha\in\Lambda}B^{\alpha})(x)=\,\{\,x\,,\,\;\bigwedge_{\alpha\in\Lambda}\mu_{V}^{\alpha}(x)\,,\;\;\bigwedge_{\alpha\in\Lambda}\nu_{V}^{\alpha}(x)\;:\;\; \forall\, x \,\in\, V\,\}$\\\\
and the \textbf{union} of those intuitionistic fuzzy hypervector spaces is defined as\\
$(\bigcup_{\alpha\in\Lambda}B^{\alpha})(x)=\,\{\,x\,,\,\;\bigvee_{\alpha\in\Lambda}\mu_{V}^{\alpha}(x)\,,\;\;\bigvee_{\alpha\in\Lambda}\nu_{V}^{\alpha}(x)\;:\;\; \forall\, x\,\in\, V\,\}$\\
\end{definition}

\begin{theorem}
The intersection of any family of intuitionistic fuzzy hypervector spaces of a hypervector space V is a intuitionistic fuzzy hypervector space.
\end{theorem}
\textbf{Proof}: Let $\{B^{\alpha}:\alpha\in\Lambda\}$ be a family of intuitionistic fuzzy hypervector spaces of V over the same intuitionistic fuzzy hyperfield $A\,=\{\,(\,x,\,\mu_{F}(x),\,\nu_{F}(x)\,)\;:\;\;x\,\in\,F\,\}$.\\
Let $B\,=\,(\bigcap_{\alpha\in\Lambda}B^{\alpha})(x)\,=\,\{\,(\,x,\,\mu_{V}(x),\,\nu_{V}(x)\,)\;:\;\;x\,\in\,V\,\}$\\
where $\mu_{V}(x)\,=\,\bigwedge_{\alpha\in\Lambda}\mu_{V}^{\alpha}(x)$ and $\nu_{V}(x)\,=\,\bigwedge_{\alpha\in\Lambda}\nu_{V}^{\alpha}(x)$\\
Let $x,\, y\,\in\, V$ and a, b $\in$ F\\
$\bigwedge_{z\in a\ast x\# b\ast y}\mu_{V}(z)=\bigwedge_{z\in a\ast x\# b\ast y}(\bigwedge_{\alpha\in\Lambda}\mu_{V}^{\alpha}(z))$\\
\smallskip$\hspace{3cm}=\bigwedge_{\alpha\in\Lambda}(\bigwedge_{z\in a\ast x\# b\ast y}\mu_{V}^{\alpha}(z)$\\
\smallskip$\hspace{3cm}\geq\bigwedge_{\alpha\in\Lambda}\{(\mu_{F}(a)\wedge\mu_{V}^{\alpha}(x))\bigwedge(\mu_{F}(b)\wedge\mu_{V}^{\alpha}(y))\}$\\
\smallskip$\hspace{3cm}=\{(\mu_{F}(a)\bigwedge(\bigwedge_{\alpha\in\Lambda}\mu_{V}^{\alpha}(x))\}\bigwedge\{\mu_{F}(b)\bigwedge
(\bigwedge_{\alpha\in\Lambda}\mu_{V}^{\alpha}(y))\}$\\
\smallskip$\hspace{3cm}=(\mu_{F}(a)\wedge\mu_{V}(x))\bigwedge(\mu_{F}(b)\wedge\mu_{V}(y))$\\
Therefore $\bigwedge_{z\in a\ast x\# b\ast y}\mu_{V}(z)\geq(\mu_{F}(a)\wedge\mu_{V}(x))\bigwedge(\mu_{F}(b)\wedge\mu_{V}(y))$\\
Again $\mu_{F}(1)\geq\mu_{V}^{\alpha}(\theta)$ , $\forall\alpha\in\Lambda$\\
Therefore $\mu_{F}(1)\geq\bigwedge_{\alpha\in\Lambda}\mu_{V}^{\alpha}(\theta)$ \\
i.e $\mu_{F}(1)\geq\mu_{V}(\theta)$\\
\textbf{Next} $\bigvee_{z\in a\ast x\# b\ast y}\nu_{V}(z)=\bigvee_{z\in a\ast x\# b\ast y}(\bigwedge_{\alpha\in\Lambda}\nu_{V}^{\alpha}(z))$\\
\smallskip$\hspace{4.2cm}\leq\,\bigwedge_{\alpha\in\Lambda}(\bigvee_{z\in a\ast x\# b\ast y}\nu_{V}^{\alpha}(z))$\\
\smallskip$\hspace{4.2cm}\leq\,\bigwedge_{\alpha\in\Lambda}\{(\nu_{F}(a)\vee\nu_{V}^{\alpha}(x))\bigvee(\nu_{F}(b)\vee\nu_{V}^{\alpha}(y))\}$\\
\smallskip$\hspace{4.2cm}=\{(\nu_{F}(a)\bigvee(\bigwedge_{\alpha\in\Lambda}\nu_{V}^{\alpha}(x))\}\bigvee\{\nu_{F}(b)\bigvee
(\bigwedge_{\alpha\in\Lambda}\nu_{V}^{\alpha}(y))\}$\\
\smallskip$\hspace{4.2cm}=(\nu_{F}(a)\vee\nu_{V}(x))\bigvee(\nu_{F}(b)\vee\nu_{V}(y))$\\
Therefore $\bigvee_{z\in a\ast x\# b\ast y}\nu_{V}(z)\leq(\nu_{F}(a)\vee\nu_{V}(x))\bigvee(\nu_{F}(b)\vee\nu_{V}(y))$\\
Again $\nu_{F}(1)\,\leq\,\nu_{V}^{\alpha}(\theta)$ , $\forall\alpha\in\Lambda$\\
Therefore $\nu_{F}(1)\,\leq\,\bigwedge_{\alpha\in\Lambda}\nu_{V}^{\alpha}(\theta)$ \\
i.e $\nu_{F}(1)\,\leq\,\nu_{V}(\theta)$\\
Therefore B is also a intuitionistic fuzzy hypervector space over A.\\
This completes the proof.\\

\begin{theorem}
The union of any family of intuitionistic fuzzy hypervector spaces of a hypervector space V is a intuitionistic fuzzy hypervector space.
\end{theorem}

\textbf{Proof}: Let $\{B^{\alpha}:\alpha\in\Lambda\}$ be a family of intuitionistic  fuzzy hypervector spaces of V over the same intuitionistic fuzzy hyperfield $A\,=\{\,(\,x,\,\mu_{F}(x),\,\nu_{F}(x)\,)\;:\;\;x\,\in\,F\,\}$.\\
Let $B\,=\,\bigvee_{\alpha\in\Lambda}\mu_{V}^{\alpha}(x)\,=\,\{\,(\,x,\,\mu_{V}(x),\,\nu_{V}(x)\,)\;:\;\;x\,\in\,V\,\} $\\
where $\mu_{V}(x)\,=\,\bigvee_{\alpha\in\Lambda}\mu_{V}^{\alpha}(x)$ and $\nu_{V}(x)\,=\,\bigvee_{\alpha\in\Lambda}\nu_{V}^{\alpha}(x)$\\
Let $x,\,y\,\in$ V and a, b $\in$ F\\
$\bigwedge_{z\in a\ast x\# b\ast y}\mu_{V}(z)=\bigwedge_{z\in a\ast x\# b\ast y}(\bigvee_{\alpha\in\Lambda}\mu_{V}^{\alpha}(z))$\\
\smallskip$\hspace{3cm}\geq\bigvee_{\alpha\in\Lambda}(\bigwedge_{z\in a\ast x\# b\ast y}\mu_{V}^{\alpha}(z))$\\
\smallskip$\hspace{3cm}\geq\bigvee_{\alpha\in\Lambda}\{(\mu_{F}(a)\wedge\mu_{V}^{\alpha}(x))\bigwedge(\mu_{F}(b)\wedge\mu_{V}^{\alpha}(y))\}$\\
\smallskip$\hspace{3cm}=\{(\mu_{F}(a)\bigwedge(\bigvee_{\alpha\in\Lambda}\mu_{V}^{\alpha}(x))\}\bigwedge\{\mu_{F}(b)\bigwedge
(\bigvee_{\alpha\in\Lambda}\mu_{V}^{\alpha}(y))\}$\\
\smallskip$\hspace{3cm}=(\mu_{F}(a)\wedge\mu_{V}(x))\bigwedge(\mu_{F}(b)\wedge\mu_{V}(y))$\\
Therefore $\bigwedge_{z\in a\ast x\# b\ast y}\mu_{V}(z)\geq(\mu_{F}(a)\wedge\mu_{V}(x))\bigwedge(\mu_{F}(b)\wedge\mu_{V}(y))$\\
Again $\mu_{F}(1)\geq\mu_{V}^{\alpha}(\theta)$ , $\forall\alpha\in\Lambda$\\
Therefore $\mu_{F}(1)\geq\bigvee_{\alpha\in\Lambda}\mu_{V}^{\alpha}(\theta)$ \\
i.e $\mu_{F}(1)\geq\mu_{V}(\theta)$\\
\textbf{Next} $\bigvee_{z\in a\ast x\# b\ast y}\nu_{V}(z)=\,\bigvee_{z\in a\ast x\# b\ast y}(\bigvee_{\alpha\in\Lambda}\nu_{V}^{\alpha}(z))$\\
\smallskip$\hspace{4.2cm}=\,\bigvee_{\alpha\in\Lambda}(\bigvee_{z\in a\ast x\# b\ast y}\nu_{V}^{\alpha}(z))$\\
\smallskip$\hspace{4.2cm}\leq\,\bigvee_{\alpha\in\Lambda}\{(\nu_{F}(a)\vee\nu_{V}^{\alpha}(x))\bigvee(\nu_{F}(b)\vee\nu_{V}^{\alpha}(y))\}$\\
\smallskip$\hspace{4.2cm}=\{(\nu_{F}(a)\bigvee(\bigvee_{\alpha\in\Lambda}\nu_{V}^{\alpha}(x))\}\bigvee\{\nu_{F}(b)\bigvee
(\bigvee_{\alpha\in\Lambda}\nu_{V}^{\alpha}(y))\}$\\
\smallskip$\hspace{4.2cm}=(\nu_{F}(a)\vee\nu_{V}(x))\bigvee(\nu_{F}(b)\vee\nu_{V}(y))$\\
Therefore $\bigvee_{z\in a\ast x\# b\ast y}\nu_{V}(z)\leq(\nu_{F}(a)\vee\nu_{V}(x))\bigvee(\nu_{F}(b)\vee\nu_{V}(y))$\\
Again $\nu_{F}(1)\,\leq\,\nu_{V}^{\alpha}(\theta)$ , $\forall\alpha\in\Lambda$\\
Therefore $\nu_{F}(1)\,\leq\,\bigvee_{\alpha\in\Lambda}\nu_{V}^{\alpha}(\theta)$ \\
i.e $\nu_{F}(1)\,\leq\,\nu_{V}(\theta)$\\
Therefore B is also a intuitionistic fuzzy hypervector space over A.\\
This completes the proof.\\

\section{Linear Transformation}

\smallskip\hspace{0.5cm}In this section we established the definition of a Linear transformation on hypervector Spaces and deduce a important theorem relative to a intuitionistic fuzzy concept.\\
\begin{definition}
Let $(V,\,\#,\,\ast)$ and $(W,\,\#^{'},\,\ast^{'})$ be two hypervector space over the same hyperfield $(F,\,\oplus,\,.)$. A mapping $T\,:V\,\rightarrow$ W is called
 Linear transformation iff\\
$(i)$ T$(x\,\#\,y)\subseteq T(x)\,\#^{'}\,T(y)$,\\
$(ii)$ T$(a\ast x)\subseteq a \ast^{'}T(x),\;\;\;\forall\, x,\, y\in$ V and a $\in$ F\\
$(iii)$ $T(\theta)\,=\,\theta^{'}$
 \end{definition}

\begin{theorem}
Let $(V,\,\#,\,\ast)$ and $(W,\,\#^{'},\,\ast^{'})$ be two hypervector space over the same hyperfield $(F,\,\oplus,\,.)$ and $T\,:V\,\rightarrow$ W be a linear transformation.
Let $B\,=\,\{\,(\,y,\;\mu_{W}(y),\;\,\nu_{W}(y)\,)\;:\;y\,\in\,W\,\}$ be a intuitionistic fuzzy hypervector space over $A\,=\,\{\,(\,a,\;\mu_{F}(a),\;\,\nu_{F}(a)\,)\;:\;a\,\in\,F\,\}$. Then $T^{-1}(B)\,=\,\{\,x,\,\;T^{-1}(\mu_{W})(x)\;,\;\,T^{-1}(\nu_{W})(x)\;:\;x\,\in\,V\,\}$ is a intuitionistic fuzzy hypervector space of V over A.
\end{theorem}

\textbf{Proof}: Let a, b $\in$ F and $\alpha, \beta\in$ V.\\
Then $\bigwedge_{x\in a\ast\alpha\#b\ast\beta}T^{-1}(\mu_{W})(x)=\bigwedge_{x\in a\ast\alpha\#b\ast\beta}\mu_{W}(T(x))$\\
\smallskip$\hspace{5.6cm}=\bigwedge_{T(x)\in T(a\ast\alpha\#b\ast\beta)}\mu_{W}(T(x))$\\
\smallskip$\hspace{5.6cm}\geq\bigwedge_{T(x)\in a\ast^{'}T(\alpha)\#^{'}b\ast^{'}T(\beta)}\mu_{W}(T(x))$\\
\smallskip$\hspace{5.6cm}\geq(\mu_{F}(a)\wedge\mu_{W}(T(\alpha))\bigwedge(\mu_{F}(b)\wedge\mu_{W}(T(\beta)))$\\
\smallskip$\hspace{5.6cm}=(\mu_{F}(a)\wedge T^{-1}(\mu_{W})(\alpha))\bigwedge(\mu_{F}(b)\wedge T^{-1}(\mu_{W})(\beta))$\\
Again $T^{-1}(\mu_{W})(\theta)=\mu_{W}(T(\theta))=\mu_{W}(\theta^{'})$ , where $\theta^{'}$ be a null vector of W.\\
\smallskip$\hspace{6cm}\leq\mu_{F}(1)$\\
i.e $\mu_{F}(1)\geq T^{-1}(\mu_{W})(\theta)$\\

And $\bigvee_{x\in a\ast\alpha\#b\ast\beta}T^{-1}(\nu_{W})(x)=\bigvee_{x\in a\ast\alpha\#b\ast\beta}\nu_{W}(T(x))$\\
\smallskip$\hspace{5.6cm}=\bigvee_{T(x)\in T(a\ast\alpha\#b\ast\beta)}\nu_{W}(T(x))$\\
\smallskip$\hspace{5.6cm}\leq\,\bigvee_{T(x)\in a\ast^{'}T(\alpha)\#^{'}b\ast^{'}T(\beta)}\nu_{W}(T(x))$\\
\smallskip$\hspace{5.6cm}\leq(\nu_{F}(a)\vee\nu_{W}(T(\alpha))\bigvee(\nu_{F}(b)\vee\nu_{W}(T(\beta)))$\\
\smallskip$\hspace{5.6cm}=(\nu_{F}(a)\vee T^{-1}(\nu_{W})(\alpha))\bigvee(\nu_{F}(b)\vee T^{-1}(\nu_{W})(\beta))$\\
Again $T^{-1}(\nu_{W})(\theta)=\nu_{W}(T(\theta))=\nu_{W}(\theta^{'})$ , where $\theta^{'}$ be a null vector of W.\\
\smallskip$\hspace{6cm}\geq\nu_{F}(1)$\\
i.e $\nu_{F}(1)\leq T^{-1}(\nu_{W})(\theta)$\\
Therefore $T^{-1}(B)$ is a intuitionistic fuzzy hypervector space of a hypervector space V over A.

\end{document}